\documentclass[11pt]{amsart}
\usepackage{amssymb, latexsym, amsmath, amsfonts, tikz}
\usepackage{hyperref, enumitem, fancyhdr}
\usepackage{color}

\newtheorem{thm}{Theorem}[section]
\newtheorem{cor}[thm]{Corollary}

\newtheorem{prop}[thm]{Proposition}
\theoremstyle{definition}
\newtheorem{defn}[thm]{Definition}
\theoremstyle{remark}

\numberwithin{equation}{section}
\theoremstyle{remark}

\newcommand{\mbb}{\mathbb}
\newcommand{\ra}{\rightarrow}

\newcommand{\sm}{\setminus}
\newcommand{\ep}{\epsilon}
\newcommand{\no}{\noindent}

\newcommand{\cal}{\mathcal}

\newcommand{\short}[1]{\textit{Short}  $\mathbb C^{#1}$}

\begin{document}
\title{Further remarks on rigidity of H\'{e}non maps}
%%\thanks{The author was supported in part by a UGC--CAS Grant}
\thanks{The author is supported by National Post-Doctoral fellowship, SERB, India}
\author{Ratna Pal}

\address{RP: Department of Mathematics, Indian Institute of Science Education and Research, Pune, Maharashtra-411008, India}
\email{ratna.math@gmail.com}
%
%\address{KV: Department of Mathematics, Indian Institute of Science, Bangalore 560 012, India}
%\email{kverma@iisc.ac.in}
%
%
\begin{abstract}
%The goal of the present article is to improve and extend the {\it{rigidity}} result of H\'{e}non maps obtained in \cite{BPK}.
 For a H\'{e}non map $H$ in $\mathbb{C}^2$, we characterize the polynomial automorphisms of $\mathbb{C}^2$ which keep any fixed level set of the Green function of $H$ completely invariant. The interior of any non-zero sublevel set of the Green function of a H\'{e}non map turns out to be a \short{2} and as a  consequence of our characterization, it follows that there exists no polynomial automorphism  apart from possibly the affine automorphisms which acts as an automorphism on any of these \short{2}'s. Further, we prove that if any two level sets of the Green functions of a pair of H\'{e}non maps coincide, then they almost commute. 
 %Let $H$ be a H\'{e}non map in $\mathbb{C}^2$.  Let  $F$ be any polynomial automorphism of degree greater than or equal to $2$ which keeps the non-escaping set of $H$, i.e., the zero level set of the Green function of $H$, completely invariant, then  either $F$ or $F^{-1}$ is a H\'{e}non map and  some iterates of $F$ and $H$ coincides.  Further, we show that there exists no polynomial automorphism of degree greater than or equal to $2$ under which any non-zero level set of the Green function of $H$ remains invariant. As a  consequence, it follows that there exists no polynomial automorphism possibly apart from the affine automorphisms which acts as an automorphism of any non-zero sublevel set of the Green function of $H$, which form an important class of \short{2}. Finally we prove that if any two level sets of the Green functions of a pair of H\'{e}non maps coincide, then they almost commute. 
\end{abstract}

\maketitle

\section{Introduction}
\subsection{Background}

\hspace{5mm}

For a polynomial $P$ in the complex plane,  the {\it Julia set} $J_P$ is the complement of the open set in $\mathbb{C}$ where the sequence $\{P^n\}_{n\geq 1}$ form a normal family locally. The dynamical behavior of the points in the Julia set is extremely chaotic  and in general the Julia set of a polynomial has a very complicated structure. It is easy to check that if two polynomials $P$ and $Q$ commute, then their Julia sets coincide. Conversely, it follows from the work of Baker--Er\"{e}menko \cite{BE} and Beardon \cite{Be1} that if $J_P=J_Q$ for two polynomials $P$ and $Q$ of degree greater than or equal to $2$, then
\[
P \circ Q = \sigma \circ Q \circ P
\]
where $\sigma(z) = az + b$ with $\vert a \vert = 1$ and $\sigma(J_P) = J_P$.  The Julia sets of two polynomial in the complex plane, of degree greater than or equal to $2$, coincide if and only if they commute upto a rigid motion in the complex plane. Now onwards, this property of Julia sets will be referred as a {\it{rigidity}} property of Julia sets. 

\medskip 
Recently, in \cite{BPK}, an analogue of the above-mentioned {\it rigidity} phenomenon has been proved for the H\'{e}non maps in $\mathbb{C}^2$  which, by the classification theorem of Freidland--Milnor (\cite{FM}),  is the most important class of automorphisms in $\mathbb{C}^2$ from the point of view of dynamics.  The class of H\'{e}non maps consists of the  polynomial automorphisms of 
$\mathbb{C}^2$ of the form
\begin{equation}\label{henon form}
H = H_m \circ H_{m-1} \circ \cdots \circ H_1
\end{equation}
where
\begin{equation*}
H_j(x, y) = (y, p_j(y) - \delta_j x) 
\end{equation*}
with $p_j$ a polynomial of degree $d_j \ge 2$ with highest degree coefficient $c_j\in \mathbb{C}$ and $0\neq \delta_j \in \mathbb{C}$. The degree of $H$ is $d = d_1d_2 \ldots d_m$.  Interested readers can look at \cite{BS1}, \cite{BS2} and \cite{BS3} for a detailed study of the dynamics of H\'{e}non maps. 

\medskip 
As in the case of polynomials in $\mathbb{C}$, one can give an analogous definition of Julia sets for  H\'{e}non maps in $\mathbb{C}^2$. For a H\'{e}non map $H$ in $\mathbb{C}^2$,  the {\it Julia sets} $J_H^\pm$  are defined as the complement of the open sets in $\mathbb{C}^2$ where the sequence $\{H^{\pm n}\}$ form normal families locally. Here $H^{\pm n}$ denote the $n$-fold iterates of $H$ and $H^{-1}$, respectively. It turns out that 
\[
J_H^\pm =\partial K_H^\pm
\]
where 
\[
K^{\pm}_H = \{(x, y) \in \mathbb C^2 : \;\text{the sequence}\; \left(H^{\pm n}(x, y) \right) \; \text{is bounded} \},
\]
the set of {\it non-escaping} points.

\begin{defn}
Let $H$ be an automorphism on $\mathbb{C}^2$. A set $S\subseteq \mathbb{C}^2$ is called completely invariant under the map $H$ if $H(S)=S$.
\end{defn}
Note that $K_H^\pm$ are completely invariant under $H$. In \cite{BPK}, for a H\'{e}non map $H$, we proved the following {\it {rigidity}} theorem. 
\subsection*{Known Theorem 1 (KT1).}
Let  $F$ be an automorphism of $\mathbb{C}^2$ which keeps the non-escaping sets $K_H^\pm$ completely invariant, then $F$ is a polynomial automorphism. If $\deg F \geq 2$, then either $F$ or $F^{-1}$ is a H\'{e}non map. Further, $F$ shares a close relation with $H$, viz., 
\begin{equation} \label{relFH1}
F^{\pm 1} \circ H = C \circ H \circ F^{\pm 1}
\end{equation}
where $C$ a linear map of the form $(x, y) \mapsto (\delta_- x, \delta_+ y)$ with $\vert \delta_\pm \vert = 1$. Also, 
\begin{equation}\label{relFH2}
\text{ either } F^m=\sigma H^n \text{ or } F^{-m}=\sigma H^n
\end{equation}
for some $m,n \in \mathbb{N}$ and for some affine automorphism $\sigma$ in $\mathbb{C}^2$.  

\medskip 
The same techniques, which are used to prove  the above theorem, gives the following version of {\it {rigidity}} theorem for Julia sets of H\'{e}non maps.

\subsection*{Known Theorem 2 (KT2).}
Let $H$ and $F$ be two H\'{e}non maps such that their Julia sets coincide, i.e, $J_H^\pm=J_F^\pm$, then (\ref{relFH1}) and (\ref{relFH2}) hold. Conversely, if $F$ and $H$ be two H\'{e}non maps satisfying 
\[
F\circ H=C\circ H \circ F
\]
where $C:(x,y)\mapsto (\delta_- x, \delta_+ y)$ with $\lvert \delta_\pm \rvert =1$ and $C(K_H^\pm)=K_H^\pm$, then $J_H^\pm=J_F^\pm$.

\medskip 
The goal of the present article is to improve and extend the {\it{rigidity}} result of H\'{e}non maps obtained in \cite{BPK}.
\subsection{Main Results}

\hspace{20mm}

\no 
Let $H$ be a H\'{e}non map. Now if we start with a polynomial automorphism $F$ of $\mathbb{C}^2$ such that $F(K_H^+)=K_H^+$, then we can recover the same relations between $H$ and $F$ as in (\ref{relFH1}) and (\ref{relFH2}) (see (KT1)). This shows that the condition $F(K_H^-)=K_H^-$ in (KT1) is redundant if we start with a polynomial automorphism $F$ in $\mathbb{C}^2$.  With these words, we present our first theorem.
\begin{thm}	\label{new rigidity}
	Let $H$ be a H\'{e}non map in $\mathbb{C}^2$ and  $F$ be a polynomial automorphism of $\mathbb{C}^2$ which keeps $K_H^+$ completely invariant.  Then, 
	\begin{itemize}
		\item[(a)]
		if $\deg(F) = 1$,  $F$ is of the form 
		\[
		(x,y)\mapsto (ax+f, dy+g)
		\]
		with $\lvert a \rvert= \lvert d \rvert=1$ and 
		\item[(b)]
		if $\deg F \geq 2$, then either $F$ or $F^{-1}$ is a H\'{e}non map and accordingly there exist $m,n \in \mathbb{N}^*$ such that
		$$
		\text{either } F^{ m}= H^n \text{ or } F^{ -m}= H^n .
		$$ 
		Further, there exists $m\in \mathbb{N}^*$ such that 
		$$
		\text{either } F^{ m}\circ H= H \circ F^m \text{ or } F^{ -m}\circ H= H\circ F^{-m} .
		$$ 
	\end{itemize}
\end{thm}
\no 
The Green functions of a H\'{e}non map $H$ is defined as follows:
\begin{equation} \label{Green}
G^{\pm}_H(x, y) := \lim_{n \rightarrow \infty} \frac{1}{d^n} \log^+ \Vert H^{\pm n}(x, y) \Vert
\end{equation}
for all $(x,y)\in \mathbb{C}^2$. Here $\log^+(x)=\max \{\log x, 0\}$. The functions $G_H^\pm$ are continuous, plurisubharmonic, non-negative on $\mbb C^2$ and pluriharmonic on $\mbb C^2 \sm K^{\pm}_H$ vanishing precisely on $K^{\pm}_H$. 

\medskip 
In case of Theorem 1.1, \cite{BPK}, since we start with an automorphism $F$ which keeps both $K_H^\pm$ invariant,  a direct analysis of the possible forms of $F$ (using Jung's theorem, \cite{J}) shows  that $F$ (or $F^{-1}$) must be a H\'{e}non map. Consequently, it follows  immediately that $G_H^+=G_F^+$ (or $G_H^+=G_F^-$). But in the present case, since only the invariance of $K_H^+$ is available, a similar analysis using Jung's theorem does not a priori gurantee  that $F$ is a H\'{e}non map, rather we get that $F$ (and hence $F^{-1}$) is a regular (hence H\'{e}non-type) map (see Section 2 for the definitions of regular maps and H\'{e}non-type maps). Then it requires some work to show that the Green functions of $H$ and $F$ coincide, i.e., $G_H^+=G_F^+$  or $G_H^+=G_F^-$ (one can define Green functions of regular and H\'{e}non-type maps in the similar fashion as in the case of H\'{e}non maps and it is discussed in Section 2).
%Theorem \ref{new rigidity} is comparable to a result due to Lamy (Theorem 5.4 in \cite{L})) where he showed that for any two H\'{e}non-type maps (the maps in $\mathbb{C}^2$ which are conjugate to  H\'{e}non maps) $F$ and $H$ if $G_H^+=G_F^+$ , then some iterates of $F$ and $H$ coincide, i.e., there exist $m,n \in \mathbb{N}$ such that $F^m=H^n$. In fact, to complete the proof of Thm \ref{new rigidity}, we use Lamy's result once we show that any polynomial automorphism  $F$ (of degree greater than or equal to $2$) satisfying  $F(K_H^+)=K_H^+$ where $H$ is a H\'{e}non map, turns out to be a H\'{e}non-type map and  their Green functions coincide, i.e., $G_H^+=G_F^+$.
Then we use Lamy's theorem (\cite{L}), to show that some iterates of $H$ and $F$ (or $F^{-1}$) coincide, i.e., there exist $m,n \in \mathbb{N}$ such that $F^m=H^n$ (or $F^{-m}=H^n$)  which shows that $F$ (or $F^{-1}$) is indeed a H\'{e}non map. In fact, in sprit Theorem \ref{new rigidity} is similar to Theorem 5.4 in \cite{L}.

\medskip 
For each $c>0$, let (see \cite{DS})
\begin{equation*}
\tilde{\Omega}_{H,c}^\pm =\left\{(x,y)\in \mathbb{C}^2: G_H^\pm(x,y)  \leq c\right\}, \:\ K_{H,c}^\pm =\left\{(x,y)\in \mathbb{C}^2: G_H^\pm(x,y)  = c\right\},\;\ J_{H,c}^\pm= \partial K_{H,c}^\pm.
\end{equation*}
Note that for $c>0$, the set $K_{H,c}^+$ has empty interior and thus 
\[
K_{H,c}^+=J_{H,c}^+.
\]
Further define
\begin{equation*}
G_{H,c}^\pm(x,y) := \max \left \{G_H^\pm(x,y) -c,0\right\}
\end{equation*}
for $(x,y)\in \mathbb{C}^2$. 
The functions $G_{H,c}^\pm$ are continuous, plurisubharmonic, non-negative on $\mbb C^2$ and pluriharmonic on $\mbb C^2 \sm J_{H,c}^\pm$ vanishing precisely on $\tilde{\Omega}_{H,c}^\pm$. 
 It can be shown that for a H\'{e}non map $H$, the {\it non-escaping} sets are the zero level sets of the Green functions $G_H^\pm$, i.e.,
$$
K_H^\pm=\{(x,y)\in \mathbb{C}^2: G_H^\pm(x,y)=0\}.
$$
Clearly, $K_{H,0}^\pm=K_H^\pm$. Theorem \ref{new rigidity} gives a characterization of automorphisms  in $\mathbb{C}^2$ in terms of $H$, which keep $K_{H,0}^\pm$ completely invariant. The next theorem shows that  in case $c>0$,  there exists no automorphism except possibly the affine automorphisms which keeps $K_{H,c}^\pm$ completely invariant.
\begin{thm} \label{Glevel}
	Let $H$ be a H\'{e}non map in $\mathbb{C}^2$. If $F$ is a polynomial  automorphism  in $\mathbb{C}^2$ such that 
	\begin{equation*}
	F\left(K_{H,c}^+\right)=K_{H,c}^+ 
	\end{equation*}
	for some $c >0$, then $F$ is an affine automorphism of the follwing form
	\[
	(x,y)\mapsto (ax+by+f, dy+g).
	\]
	%\end{itemize}
\end{thm}
 In \cite{F}, Forn{\ae}ss showed the existence of so called \short{k}. A  domain $\Omega$ which can be expressed as an increasing union of unit balls (upto biholomorphism) such that Kobayashi metric vanises identically in $\Omega$,  but allows a bounded (above) pluri-subharmonic function, is called  \short{k}.  For a H\'{e}non map $H$, it can be shown that the interior of any non-zero sublevel set of the Green function $G_H^+$, i.e.,
\[
\Omega_{H,c}=\left\{(x,y)\in \mathbb{C}^2: G_H^+(x,y)  <  c\right\}
\]
is a  \short{2}, for any $c>0$ (see \cite{F}). Since $\Omega_{H,c}$ is essentially an increasing union of Euclidean balls in $\mathbb{C}^2$ whose automorphism group is well-understood,  it is an interesting task to understand the automorphism group of $\Omega_{H,c}$. Since any polynomial automorphism $F$ in $\mathbb{C}^2$ which acts as an automorphism of $\Omega_{H,c}$ will keep $K_{H,c}^+$ completely invariant,  a simple application of Theorem \ref{Glevel} gives the following proposition. 
\begin{prop}
	For any $c>0$, there exists no polynomial automorphism of $\mathbb{C}^2$ except possibly the affine automorphisms 
	of the form
	\[
	(x,y)\mapsto (ax+by+f, dy+g).
	\] 
	which acts as an automorphism of $\Omega_{H,c}$.
\end{prop}
It follows from Theorem 1.1 in \cite{BPK} that if the zero level sets of Green functions of two H\'{e}non maps (or the Julia sets) coincide, then they almost commute. We prove  that the same is true if any two level sets of the Green functions of a pair of H\'{e}non maps coincide.  	

\begin{thm} \label{two levels}
	Let $H_1$ and $H_2$ be two H\'{e}non maps of degree $d_{H_1}$ and $d_{H_2}$, respectively such that 
	\[
	J_{H_1,c_1}^+= J_{H_2,c_2}^+ \text{ and } J_{H_1,d_1}^-= J_{H_2,d_2}^-
	\]	
	for some $c_1, c_2, d_1, d_2 \geq 0$, then 
	$$
	H_2\circ H_1=C \circ H_2 \circ H_1.
	$$
	Here $C(x,y)=(\delta_+ x, \delta_- y)$ with $\lvert \delta_+\rvert= e^ {c(d_{H_1}-1) (d_{H_2}-1)}$ and  $\lvert \delta_-\rvert=e^ {d(d_{H_1}-1) (d_{H_2}-1)}$ where $c=c_1-c_2$ and $d=d_1-d_2$.
\end{thm}
 Before we start proving our main theorems, we gather a few preparatory stuff in the next section.  Proofs of Theorem \ref{new rigidity}, Theorem \ref{Glevel} and Theorem \ref{two levels} appear in Section 3, Section 4 and Section 5, respectively.

\section{Preliminaries}	
Readers are referred to \cite{BS1}, \cite{BS2}, \cite{BS3} and \cite{DS} for a detailed study of the material in this section.

\medskip 
 For $R>0$, let us first define a filtration of $\mathbb{C}^2$ as follows:
 \begin{align*}
 V^+_R &= \{ (x,y) \in \mathbb C^2: \vert x \vert < \vert y \vert, \vert y \vert > R \},\\
 V^-_R &= \{ (x,y) \in \mathbb C^2: \vert y \vert < \vert x \vert, \vert x \vert > R \},\\
 V_R &= \{ (x, y) \in \mathbb C^2: \vert x \vert, \vert y \vert \le R \}.
 \end{align*}
 For a given H\'{e}non map $H$ of degree $d$, there exists $R > 0$ such that
 \[
 H(V^+_R) \subset V^+_R, \; H(V^+_R \cup V_R) \subset V^+_R \cup V_R,
 \]
  \[
 H^{-1}(V^-_R) \subset V^-_R, \; H^{-1}(V^-_R \cup V_R) \subset V^-_R \cup V_R,
 \]
 \[
   K_H^{\pm} \subset V_R \cup V^{\mp}_R  \text{ and } \mathbb C^2 \sm K^{\pm}_H = \bigcup_{n=0}^{\infty} (H^{\mp n})(V^{\pm}_R).
  \]
  Recall that 
  \[
  K^{\pm}_H = \{(x, y) \in \mathbb C^2 : \;\text{the sequence}\; \left(H^{\pm n}(x, y) \right) \; \text{is bounded} \}.
  \]
  As defined in the previous section, the Green functions
 \[
 G^{\pm}_H(x, y) = \lim_{n \rightarrow \infty} \frac{1}{d^n} \log^+ \Vert H^{\pm n}(x, y) \Vert
 \]
 for $(x,y)\in \mathbb{C}^2$. The functions $G_H^\pm$ are continuous, plurisubharmonic, non-negative on $\mbb C^2$ and pluriharmonic on $\mbb C^2 \sm K^{\pm}_H$ vanishing precisely on $K^{\pm}_H$. By construction, the following functorial property holds:
 \[
 G^{\pm}_H \circ H = d^{\pm 1} G^{\pm}_H.
 \]
 Both $G^{\pm}_H$ have logarithmic growth near infinity, i.e., there exists $R>0$ such that   
 \begin{equation}\label{L1}
 G_H^+ (x,y)= \log^+ \lvert y \rvert+ O(1)
 \end{equation}
 in  $\overline{V_R^+ \cup V_R}$, and 
 \begin{equation}\label{L2}
 G_H^- (x,y)= \log^+ \lvert x \rvert+ O(1)
 \end{equation}
 in $\overline{V_R^- \cup V_R}$. Hence
 \begin{equation}\label{L3}
 G_H^\pm (x,y)\leq \max \{\log^+ \lvert x\rvert, \log^+ \lvert y \rvert \}+ C
 \end{equation}
 for all $(x,y)\in \mathbb{C}^2$ and for some $C>0$. 
 %It turns out that $G_H^\pm$ are the pluricomplex Green's functions for $K^{\pm}_H$ respectively. 
 The supports of the positive closed $(1,1)$ currents 
 \[
 \mu^{\pm}_H = dd^c G^{\pm}_H
 \]
 are $J^{\pm}_H$ and $\mu_H = \mu^+_H \wedge \mu^-_H$ is an invariant measure for $H$.  
 
 \medskip 
 Recall that, for each $c>0$, we define (see \cite{DS})
 \begin{equation*}
 \tilde{\Omega}_{H,c}^\pm =\left\{(x,y)\in \mathbb{C}^2: G_H^\pm(x,y)  \leq c\right\}, K_{H,c}^\pm =\left\{(x,y)\in \mathbb{C}^2: G_H^\pm(x,y)  = c\right\} , J_{H,c}^\pm= \partial K_{H,c}^\pm
 \end{equation*}
 and recall that $K_{H,c}^+=J_{H,c}^+$ for $c>0$. Further define
 \begin{equation*}
 G_{H,c}^\pm(x,y) = \max \left \{G_H^\pm(x,y) -c,0\right\}
 \end{equation*}
 for $(x,y)\in \mathbb{C}^2$. 
 The functions $G_{H,c}^\pm$ are continuous, plurisubharmonic, non-negative on $\mbb C^2$ and pluriharmonic on $\mbb C^2 \sm J_{H,c}^\pm$ vanishing precisely on $\tilde{\Omega}_{H,c}^\pm$. Clearly, $G_{H,c}^\pm$ satisfy the same inequalities as in (\ref{L1}), (\ref{L2}) and (\ref{L3}). Further, The supports of the positive closed $(1,1)$ currents 
 \[
 \mu_{H,c}^{\pm} = dd^c G_{H,c}^{\pm}
 \]
 are $J^{\pm}_{H,c}$.
 
 \medskip
The following  theorem proved by Dinh--Sibony (see \cite{DS}) have been crucially used to establish the main theorems of the present article. 
 \begin{thm}\label{DS1}
 The current $\mu_H^+$ is the unique closed positive $(1,1)$ current of mass $1$ supported on $J_H^+$.  For any $c>0$, the current $\mu_{H,c}^+$ is a closed  positive $(1,1)$ current of mass $1$ supported on $J_{H,c}^+$.
 \end{thm}
 Any H\'{e}non map extends meromorphically to $\mathbb{P}^2$ with an isolated indeterminacy point at $I^+ = [1:0:0]$ and similarly, $H^{-1}$ extends to $\mathbb{P}^2$ with a lone indeterminacy point at $I^{-} = [0:1:0]$.
 The class of H\'{e}non maps form the most important class of {\it regular} maps in $\mathbb{C}^2$.  
 
 \medskip 
 For a polynomial $f$ in $\mathbb{C}^2$, let $\hat{f}$ and $\hat{f}^{-1}$ be the meromorphic extensions of $f$ and $f^{-1}$ to $\mathbb{P}^2$, respectively. Let $I_f^+$ and $I_f^-$ be the inderminacy points of $\hat{f}$ and 
$\hat{f}^{-1}$ in $\mathbb{P}^2$, respectively.
 \begin{defn}
 	We say that $f$ is regular if $I_f^+ \cap I_f^{-}=\emptyset$.
 \end{defn}
For a regular map $f$ in $\mathbb{C}^2$ of degree $d$,  the Green functions
\[
G^{\pm}_f(x, y) := \lim_{n \rightarrow \infty} \frac{1}{d^n} \log^+ \Vert f^{\pm n}(x, y) \Vert
\]
for $(x,y)\in \mathbb{C}^2$.  We define
 \[
K^{\pm}_f = \{(x, y) \in \mathbb C^2 : \;\text{the sequence}\; \left(f^{\pm n}(x, y) \right) \; \text{is bounded} \} \text{ and } J_f^{\pm}=\partial K_f^\pm.
\]
The functions $G_f^\pm$ are continuous, plurisubharmonic, non-negative on $\mbb C^2$ and pluriharmonic on $\mbb C^2 \sm K^{\pm}_f$ vanishing precisely on $K^{\pm}_f$. By construction, the following functorial property holds:
\[
G^{\pm}_f\circ f = d^{\pm 1} G^{\pm}_f
\]
where $d$ is the degree of $f$. Further, the functions $G^{\pm}_f$ have logarithmic growth near $I_f^{\mp}$ and the similar inequalities as in (\ref{L1}), \ref{L2} and (\ref{L3}) hold for $f$. See \cite{SW} (Section 2) for the following proposition.
\begin{prop}\label{relegular attraction}
The points $I_f^+$ and $I_f^-$ are the attracting fixed points for $f^{-1}$ and $f$, respectively. Futhermore, for any point $z\in \mathbb{C}^2\setminus K_f^\pm$,
\[
f^{\pm n} (z) \rightarrow I_f^\mp
\]	
as $n\rightarrow \infty$.
\end{prop}
The supports of the positive closed $(1,1)$ currents 
\[
\mu^{\pm}_f = dd^c G^{\pm}_f
\]
are $J^{\pm}_f$. 

\medskip 
The following theorem is due to Dinh--Sibony (\cite{DS}) 
\begin{thm}\label{DS2}
	The current $\mu_f^+$ is the unique closed $(1,1)$ current of mass $1$ supported on the sets $K_f^+$ and $J_f^+$. 
\end{thm}

\begin{defn}
	A polynomial automorphism $f$ in $\mathbb{C}^2$ is called H\'{e}non-type if 
	\[
	f=\varphi \circ h \circ \varphi^{-1}
	\]
	where $h$ is  composition of H\'{e}non maps and $\varphi$ is a polynomial automorphism in $\mathbb{C}^2$.
\end{defn}
Clearly, a regular map is a H\'{e}non-type map.

\medskip 
Let $\mathcal{PSH}(\mathbb{C}^2)$ be the collection of all pluri-subharmonic functions in $\mathbb{C}^2$. Set 
\[
\mathcal{L}=\{v\in \mathcal{PSH}(\mathbb{C}^2): v(z)\leq \log^+ \lVert z \rVert+M \text{ with some } M>0\}
\]
where $\log^+(x)=\max \{\log x, 0\}$.
\begin{defn}
	For a subset $S\subseteq \mathbb{C}^2$, the function 
	\[
	L_S(z):=\sup \{u(z): u\in \mathcal{L}, u\leq 0 \text{ on } S\}
	\]
	for $z\in \mathbb{C}^2$, is called the pluricomplex Green function of $S$.
\end{defn}
One can look at Proposition 8.4.10 in \cite{MNTU} for the proof of the following proposition. 
\begin{prop}\label{pluri henon}
	The pluricomplex Green functions of the sets $K_H^\pm$ and of the sets $J_H^\pm$ are $G_H^\pm$.
\end{prop}

\section{Proof of Theorem \ref{new rigidity}}
 Suppose that $\deg(F) \ge 2$.  Then we show that
	\begin{center}
		{\it $F$ is a regular automorphism:}
	\end{center}
	
	%\subsection*{Case 1:} Suppose that $\deg(F) \ge 2$
	
	%\medskip 
	\no 
	That $F$ is a regular automorphism is obtained  following the same line of arguments as in the  proof of Theorem 1.1 in \cite{BPK}, which shows that any polynomial automorphism preserving the non-escaping sets $K_H^\pm$  is essentially a H\'{e}non map.  In present case, due to unavilability of invariance of $K_H^-$ under $F$,  we need to strech the arguments given in \cite{BPK} accordingly to conclude that $F$ is a regular polynomial automorphism. 
	%Later using the main result in \cite{BF}, we shall be able to conclude that $F$ (or $F^{-1}$) is indeed a H\'{e}non map.   
	
	\medskip 
   By Jung's theorem (see \cite{J}), $F$ can be written as a composition of affine maps and elementary maps in $\mathbb{C}^2$. Recall that an elementary map is of the form
	\[
	e(x, y) = (\alpha x + p(y), \beta y + \gamma)
	\]
	where $\alpha \beta \not= 0$ and $p(y)$ is a polynomial in $y$. We consider following cases.
	
	\medskip
	\no
	{\it{Case (i)}}: Let
	\[
	F=a_1\circ e_1 \circ a_2 \circ e_2\circ \cdots \circ a_k \circ e_k
	\]
	for some $k\geq 1$ where the $a_i$'s are non-elementary affine maps and the 
	$e_i$'s are non-affine elementary maps.
	Without loss of generality, suppose that
	\[
	F=a_1\circ e_1 \circ a_2 \circ e_2.
	\]
	Let 
	\[
	a_1(x,y)=(\alpha_1 x+ \beta_1 y+ \delta_1,\alpha_2 x+ \beta_2 y+ \delta_2 ).
	\]
	for $\alpha_2\neq 0$. Now consider the maps
	\[
	a_1^2(x,y)=(\alpha_2 x+ \beta_2 y+ \delta_2, s_2y+r_2)
	\]
	and
	\begin{equation}\label{aff}
	a_1^1(x,y)=(bx+cy,y)
	\end{equation}
	where $b \not= 0, c=\alpha_1/\alpha_2$, $r_2=(\delta_1-c\delta_2)/b$ and $s_2=(\beta_1-c\beta_2)/b$. Note that $a_1= a_1^1\circ \tau\circ a_1^2$ where $\tau(x,y)=(y,x)$ for any $b\neq 0$. Expressing  $a_2$ in a similar way, it follows that
	\[
	F=a_1^1\tau a_1^2 e_1 a_2^1\tau a_2^2 e_2.
	\] 
	Now $a_1^2 e_1 a_2^1$ and $a_2^2 e_2$ are elementary maps, $E_1$ and $E_2$ respectively, say. Therefore,
	\[
	F=a_1^1 \tau E_1 \tau E_2
	\]
	where $E_i(x,y)=(m_i x+p_i(y), k_i y+ r_i)$ and the $p_i$'s are polynomials in 
	$y$ of degree at least $2$ for $i=1,2$.  Since $\tau E_1$ and $\tau E_2$ are H\'{e}non maps, it follows that $[1:0:0]$ is an indeterminacy point of $F$. But $F^{-1}([1:0:0])=[1:0:0]$. Thus, in this case, $F$ is regular map with $[1:0:0]$ as the forward indeterminacy point. Note that the point $[w:1:0]$ is the indeterminacy point of $F^{-1}$ for some $w\in \mathbb{C}$.
	
	\medskip
	\no 
	{\it{Case (ii):}} Let 
	\[
	F=a_1\circ e_1\circ a_2 \circ\cdots\circ e_{k-1} \circ a_k
	\]
	for some $k\geq 2$. That $F$ can not be of this form, provided $F(K_H^+)=K_H^+$, follows exactly the same set of arguments as in the proof of Theorem 1.1 in \cite{BPK} (or case (ii) in Theorem \ref{Glevel} in the present paper).
	
	\medskip
	\no 
	{\it{Case (iii):}} Let 
	\[
	F=e_1\circ a_1 \circ e_2 \circ a_2 \circ \cdots \circ e_k \circ a_k
	\]
	for some $k\geq 1$. Note that $F^{-1}$ has a form as in Case (i). Since $F^{-1}$ also keeps $K_H^+$ invariant, it follows that $F^{-1}$ is a regular map with $[1:0:0]$ as the indeterminacy point. Hence, $F$ is a regular map
	
	\medskip
	\no 
	{\it{Case (iv):}} Let 
	\[
	F=e_1\circ a_1 \circ e_2 \circ a_2 \circ \cdots\circ a_{k-1}\circ e_{k-1} \circ e_k 
	\]
	for some $k\geq 1$. For simplicity, we work with 
	\[
	F=e_1\circ a_1 \circ e_2
	\]
	and as in the previous cases, we  write
	\[
	F=e_1 a_1^1\tau a_1^2 e_2
	\]
	and thus,
	\[
	\tau F=\tau e_1 a_1^1\tau a_1^2 e_2.
	\]
	Note that both $\tau e_1 a_1^1$ and $\tau a_1^2 e_2$ are H\'{e}non maps.  Thus $\tau F$ is a H\'{e}non map. Therefore,  $F[0:1:0]=[1:0:0]$ and consequently, $F(V_R^+)$ will intersects $K_H^+$ since $[1:0:0]$ is the limit point of $K_H^+$ in $\mathbb{P}^2$. This implies, $V_R^+ \cap K_H^+\neq \emptyset$ since $F(K_H^+)=K_H^+$ . This is clearly a contradiction. Therefore, $F$ can not be of this form. 
	
	\medskip 
	Thus we prove that $F$ is a regular map. Further, the indeterminacy point of $F$ is either $[1:0:0]$ or $[w:1:0]$ and accordingly the indeterminacy point of $F^{-1}$ is either $[w:1:0]$ or $[1:0:0]$.
	\begin{center}
	{\it Green functions of $H$ and $F$ coincide:} 
	\end{center}
	 Note that in the previous section, we have shown that if $F(K_H^+)=K_H^+$, then the forms appeared in Case (i) and Case (iii) are the two possible forms of $F$. Now in Case (i), $I_F^+=[1:0:0]$ is the indeterminacy point of $F$. Hence $I_F^+$ is the attracting fixed point for $F^{-1}$ (see  Proposition \ref{relegular attraction}, Section 2).
	% Now following \cite{SW}, I think one can show that there exists $R$ sufficiently large such that
	%\[
	%K_F^+ \cap V_R^+ =\emptyset.
	%\]
	Now we have $F(K_H^+)=K_H^+$. If $z\notin K_F^+$, then $F^n(z)\rightarrow I_F^-$ as $n\rightarrow \infty$ where $I_F^{-}=[w:1:0]$ is the indeterminacy point of $F^{-1}$. But $\overline { K_H^+}=K_H^+ \cup I^+$ where $I^+=[1:0:0]$. Therefore, $K_H^+ \subset K_F^+$. Using Dinh--Sibony {\it rigidity} result (see Theorem \ref{DS2}, Section 2) for regular maps, we can conclude that 
	\[
	J^+=J_H^+=J_F^+.
	\]
	Since $H$ is a H\'{e}non map, $G_H^+$ is the pluricomplex Green function of $J_H^+$ (see \ref{pluri henon}, Section2).  Further we claim that
	 {\it{$G_F^+$ is the pluricomplex Green function of $J_F^+$.}} For $\epsilon>0$ sufficiently small, let
	 \[
	 U_F^-=\left\{(x,y)\in \mathbb{C}^2: \lvert x \rvert < (\lvert c \rvert+ \epsilon )\lvert y \rvert , \lvert y \rvert >R   \right\} 
	 \] 
	 which is clearly away from the point $[1:0:0]$. 
	 By Theorem 8.4 in \cite{DS},  
	\begin{equation*}
	\log \lVert z \rVert-M_2 \leq G_F^+(z)\leq \log \lVert z \rVert+M_1
	\end{equation*}
	for some $M_1, M_2>0$ and for all $z\in U_F^-$.  
	Fix $x_0\in \mathbb{C}$ and consider the complex line $L_{x_0}=\{(x_0,y): y\in \mathbb{C}\}$. For $R>0$ sufficiently large 
	\begin{equation*}
	\log \lVert z \rVert-K_2 \leq G_H^+(z)\leq \log \lVert z \rVert+K_1
	\end{equation*}
	in $V_R^+$ (see (\ref{L1})) for some $K_1, K_2 >0$. Now note that 
	\begin{equation*}
	L_{x_0}^R=\{(x_0,y): y\in \mathbb{C}, \lvert y \rvert>R\} \subseteq U_F^- \cap V_R^+
	\end{equation*}
	and therefore $G_H^+-G_F^+$ is bounded at infinity along the line $L_{x_0}$. Further, since $J^+=J_H^+=J_F^+$ and the function $G_H^+-G_F^+$ is harmonic in $L_{x_0}\setminus J^+$ which vanishes identically on $J^+$, it follows that
	\begin{equation*}
	G_H^+ \leq G_F^+.
	\end{equation*}
	Since 
	\[
G_F^+(z)\leq \log \lVert z \rVert+M
	\]
	for some $M>0$ and for all $z\in \mathbb{C}^2$ and $G_F^+$ vanishes identically on $J_F^+$, 
	\begin{equation*}
	G_F^+ \leq G_H^+
	\end{equation*} 
	in $\mathbb{C}^2$.
	Therefore, the Green function of $H$ and the Green function of $F$ coincide, i.e., 
	$$
	G_H^+=G_F^+
	$$
	in $\mathbb{C}^2$.
	
	\medskip 
	In the other case, that is, if $[1:0:0]$ is the indeterminacy point of $F^{-1}$, using the similar set of arguments as before, it follows that
	\[
	G_H^+=G_{F^{-1}}^+.
	\]
	
	\begin{center} 
		{\it Some iterates of $F$ and $H$ agree:}
	\end{center}
Note that since $F^{\pm 1}$ are regular maps,  they are H\'{e}non-type maps, i.e., $F^{\pm 1} $ are conjugate to some H\'{e}non maps. Further, without loss of generality, we assume that $G_H^+=G_F^+$. Therefore, by Theorem. 5.4 in \cite{L}, it follows that there exists $m,n\in \mathbb{N}^*$ such that 
\begin{equation}\label{eq iterates}
F^m=H^n.
\end{equation}

\medskip 	
\begin{center} 
	{\it $F$ is a H\'{e}non map:  }
\end{center}
Since root of a H\'{e}non map is also a H\'{e}non map (see Theorem. 4.1, \cite{BF}), it follows from (\ref{eq iterates}) that $F$ is a H\'{e}non map.

\medskip 
\no 
\begin{center} 
	{\it Some iterates of $F$ commutes with $H$:}
\end{center}
\no 
It follows from (\ref{eq iterates}) that
$$
\text{either } F^{ m}\circ H= H \circ F^m \text{ or } F^{ -m}\circ H= H\circ F^{-m}
$$ 
for some $m\in \mathbb{N}^*$.

	\begin{center} 
	{\it Description of linear automorphisms which keeps $K_H^+$ invariant:}
\end{center}
\no 	
Let $\sigma$ be an affine automorphism of the form
\[
\sigma(x,y)=(ax+by+f,cx+dy+g)
\]
such that $\sigma(K_H^+)=K_H^+$.
Thus if we take a sequence $\{(x_n,y_n)\}_{n\geq 1}\subseteq  K_H^+$ which converges to $[1:0:0]\in \mathbb{P}^2$, then
\[
(a x_n+b y_n+ f, c x_n +d y_n +g) \rightarrow [1:0:0]
\]
as $n\rightarrow \infty$, which in turn gives that $c=0$. Hence 
\[
\sigma(x,y)=(ax+by+f,dy+g).
\]
Now since $\sigma \circ H (K_H^+)=K_H^+$ and $\deg (\sigma \circ H) \geq 2$, it follows from previous description that $\sigma \circ H$ is  H\'{e}non and thus $b=0$. Therefore 
\[
\sigma(x,y)=(ax+f,dy+g).
\]
Let $(x,y)\in K_H^+$, then 
\begin{equation} \label{sigma n}
\sigma^n(x,y)=\left (a^n x+f \frac{(a^n-1)}{a-1}, d^n y + g \frac{(d^n-1)}{(d-1)}\right).
\end{equation}
Now note that if $\lvert a \rvert \leq 1$, then $\lvert d \rvert \leq 1$ since $\sigma^n (K_H^+)=K_H^+\subseteq V_R \cup V_R^-$ for all $n\geq 1$. Choose $(x_n,y_n)\in K_H^+$ such that $\lvert x_n \rvert \ra \infty$  and ${y_n}/{x_n} \ra 0$ as $n\ra \infty$.
As in (\ref{sigma n}), we have 
\begin{equation} \label{sigma nn}
\sigma^n(x_n,y_n)=\left (a^n x_n+f \frac{(a^n-1)}{a-1}, d^n y_n + g \frac{(d^n-1)}{(d-1)}\right).
\end{equation}
If $\lvert a \rvert <1$, then it follows from (\ref{sigma nn}) that $\sigma^n (x_n,y_n)\ra [0:0:0]$ as $n\ra \infty$ which is a contradiction since $\overline{ K_H^+}=K_H^+ \cup I^+$. Thus $\lvert a \rvert\geq 1$. Now 
\[
\sigma^{-1}(x,y)=\left (\frac{x}{a}-\frac{f}{a},\frac{y}{d}-\frac{g}{d} \right).
\]
Since $\sigma(K_H^+)=K_H^+$, applying same argument as before $\lvert a^{-1}\rvert \geq 1$.  Thus we get $\lvert a\rvert=1$. We have already proved that if $\lvert a\rvert \leq 1$, then $\lvert d \rvert \leq 1$ and thus we have $\lvert d \rvert=1$. So finally we get that 
\[
\sigma(x,y)=(ax+f,dy+g)
\]
with $\lvert a \rvert=\lvert d \rvert=1$.

\begin{cor}
Let $H$ and $F$ be two H\'{e}non maps such that $J_H^+=J_F^+$, then there exist $m,n \in \mathbb{N}$ such that $F^m=H^n$.
	\end{cor}

\section {Proof of the Theorem \ref{Glevel}}
Before starting the proof of \ref{Glevel}, we state the following propostion which we shall require later. The proof of the proposition follows exactly as the proof Proposition \ref{pluri henon}, hence we omit the proof.
\begin{prop} \label{pluri-complex}
	For any $c>0$, the functions $G_{H,c}^\pm$ are the pluricomplex Green functions of the sets $\tilde{\Omega}_{H,c}^\pm$ and  of the sets $K_{H,c}^\pm=J_{H,c}^\pm$.	
\end{prop} 

\subsection*{Proof of the Theorem \ref{Glevel}:} 
 Let $F$ be a polynomial automorphism of $\mathbb{C}^2$ such that $F(K_{H,c}^+)=K_{H,c}^+$, for some $c>0$. Then we prove the following equalities: 
\begin{equation}\label{Gc1}
G_{H,c}^+ \circ F^{\pm 1} (x,y)= b^\pm G_{H,c}^+(x,y)
\end{equation}
	for $(x,y)\in \mathbb{C}^2$.
	%and 
	%\begin{equation}\label{Gc2}
%	G_{H,c_2}^- \circ F^{\pm 1} (x,y)= d^\pm G_{H,c_2}^-(x,y)
%	\end{equation} 
	%for $(x,y)\in \mathbb{C}^2$. 
	Note that if (\ref{Gc1}) holds, then $b^-={(b^+)}^{-1}$ with $b^+>0$. The idea of the proof is due to Buzzard--Forn{\ae}ss (\cite{BF}).
	
	\medskip 
	 Since $F(K_{H,c}^+)=K_{H,c}^+$, it follows that $G_{H,c}^+\circ F=0$ on $K_{H,c}^+$. For any $x \in \mathbb{C}$ and consider 
	\[
	g_{x}(y):=G_{H,c}^+\circ F (x,y)
	\]
	for $y\in \mathbb{C}$. The function $g_x$ is harmonic outside the compact set $K_{H,c}^+ \cap (\{x\}\times \mathbb{C})$ and thus it is harmonic outside a large disk of radius $R>0$. Let $h_{x}$ be the harmonic conjugate of  $g_{x}$  in $\{\vert y \vert > R\}$ with period $p_{x}$. Therefore
	\[
	\psi_{x}(y)=g_{x}(y)-p_{x}\log \lvert y\rvert+ i h_{x}(y)
	\]
	is holomorphic in $\{\lvert y\rvert >R\}$. Since
	\[
	\left\lvert\exp(-\psi_{x}(y))\right\rvert \leq {\lvert y\rvert}^{p_{x}},
	\]
	the function $\exp(-\psi_{x}(y))$ has at most a pole at infinity and thus,
	\[
	\exp(-\psi_{x}(y))=y^k \exp f(y)
	\]
	where $f$ is a holomorphic function in $\{\lvert y \rvert >R\}$ having a removable singularity at infinity. 
	%Here $k = k_x$ is a positive integer. 
	Taking absolute values and then log, we get the following:
	\[
	g_{x}(y)-p_{x}\log \lvert y \rvert=-k \log \lvert y \rvert-{\rm{Re}}(f(y))
	\] 
	in $\{\lvert y\rvert >R\}$. Therefore, 
	\[
	g_{x}(y)=b_{x}\log \lvert y\rvert +O(1) 
	\]
	in $\{\lvert y\rvert >R\}$ and since $g \ge 0$ in $\mathbb{C}$,  
	\[
	g_{x}(y)=b_{x}\log^+ \lvert y\rvert +O(1) 
	\]
	in $\mathbb{C}$. 
	%Since the $O(1)$ term is bounded near infinity, $b_x$ is in fact 
	%the period of $g_x(y)$. 
	
	\medskip
	 We prove that $b_x$ is independent of $x$. To prove this, we work in a small neighbourhood of a fixed $x_0$ and  consider $R>0$  large enough. Let $p, q$ be two distinct points near $x_0$ and let $I$ be the straight line segment joining them. Then 
	\[
	\Sigma = \{(x, y) : x \in I, \vert y \vert = R \}
	\]
	is a smooth real 2-surface with two boundary components namely, 
	\[
	\{(p, y) : \vert y \vert = R \} \cup \{ (q, y): \vert y \vert = R \}. 
	\]
	We get 
	\[
	b_p - b_q = \int_{\partial \Sigma} d^c (G^+_{H,c} \circ F) = \int_{\Sigma} dd^c(G^+_{H,c} \circ F) = 0
	\]
applying Stokes' theorem. 	Here the last equality holds due to the pluriharmonicity of $G^+_{H,c} \circ F$ on $V^+_R$. Thus $b_x$ is locally constant and therefore constant everywhere. Let us write $b_x = b^+$ for all $x \in \mathbb{C}$.
	
	\medskip 
	 Now  we have 
	$
	g_{x_0}(y)=b^+ \log^+ \lvert y \rvert+ O(1) \text{ and } G_{H,c}^+(x,y)=\log^+ \lvert y\rvert +O(1)
	$
	 in $V_R^+$. Further, the difference $g_{x_0}(y)-b^+ G_{H,c}^+(y)$ is  harmonic at each $y$ for which $(x_0,y)\in \mathbb{C}^2\setminus K_{H,c}^+$ with a removable singularity at infinity and vanishes for $(x_0,y)\in K_{H,c}^+$. Therefore, $g_{x_0}(y)=b^+G_{H,c}^+(x_0,y)$ for each $y\in \mathbb{C}$. 
	Applying the same argument we get that $G_{H,c}^+\circ F-b^+G_{H,c}^+ \equiv 0$ in $\Delta(x_0;r_0)\times \mathbb{C}$. Since the difference is pluriharmonic in $\mathbb{C}^2\setminus K_{H,c}^+$  and it vanishes in $K_{H,c}^+$, we have $G_{H,c}^+\circ F=b^+G_{H,c}^+$ in $\mathbb{C}^2$. Using similar arguments we get that $$G_{H,c}^+\circ F^{-1}=b^-G_{H,c}^+$$ in $\mathbb{C}^2$ where $b^-={(b^+)}^{-1}$.

	\medskip 
	 \no 
	 Since for any $c>0$, $\overline{K_{H,c}^\pm}=K_{H,c}^\pm \cup I^\pm$ in $\mathbb{P}^2$ (see \cite{DS}), if $\deg F=1$, using the similar arguments as in Theorem \ref{new rigidity}, we get that 
	 \[
	 F(x,y)=(ax+by+f,dy+g)
	 \]
	 for $(x,y)\in \mathbb{C}^2$. In case $\deg F\geq 2$, we prove that
	\begin{center}
		{\it{ $F$ is a regular automorphism:}}
	\end{center}

	%Let $F$ be a polynomial automorphism in $\mathbb{C}^2$ such that 
	%\[
	%F(K_c^+)=K_c^+.
	%\]	
	%{\it{Claim:}} $F$ is a regular polynomial automorphism.
	
	%\medskip 
	To prove that $F$ is a regular polynomial automorphism, we shall use the similar set of arguments as in the first part of the proof of Theorem \ref{new rigidity}.  As before, the following cases arise:
	
	\medskip 
	\no 
	{\it{Case (i):}}  
	Let
	\[
	F=a_1\circ e_1 \circ a_2 \circ e_2\circ \cdots \circ a_k \circ e_k
	\]
	for some $k\geq 1$ where the $a_i$'s are non-elementary affine maps and the 
	$e_i$'s are non-affine elementary maps. As it is shown in Theorem {\ref{new rigidity}},  in this case , $F$ is a regular map. 
	
	\medskip 
	\no 
	{\it{Case (ii)}}  Let 
	\[
	F=a_1\circ e_1\circ a_2 \circ\cdots\circ e_{k-1} \circ a_k
	\]
	for some $k\geq 2$.
	For simplicity, assume that $F=a_1\circ e_1 \circ a_2$. As in the previous case we can write
	\[
	F=a_1^1 \tau a_1^2 e_1 a_2^1 \tau a_2^2
	\] 
	where $a_1^1(x,y)=(bx+cy,y)$ and $ \tau a_2^2(x,y)=(s_2 y+r_2, \alpha_2 x+ \beta_2 y+ \delta_2)$. 
	That $F$ can not be of this form follows exactly as in the proof of Theorem 1.1 in \cite{BPK}. However, since in our present case $K_{H,c}^+$ is invariant under $F$ instead of the invariance of the non-escaping set $K_H^+$ (as in Theorem 1.1 in \cite{BPK}), we need to modify our proof accordingly.  
	
	\medskip 
	 Note that  for any given $c>0$, there exists $R_c>0$ sufficiently large such that 
	\begin{equation} \label{JHc}
	K_{H,c}^+ \cap V_{R_c}^+=\emptyset. 
	\end{equation} 
	Let for each $n\in \mathbb{N}$, there exists $(x_n,y_n)\in K_{H,c}^+ \cap V_n^+$. Now by Lemma 6.3 in \cite{DS}, $\overline{ K_{H,c}^+}=K_{H,c} \cup I^+$ in $\mathbb{P}^2$. Therefore $[x_n: y_n :1]\rightarrow [1:0:0]$ as $n\rightarrow \infty$ which contradicts the fact that $(x_n,y_n)\in V_n^+$ for each $n\geq 1$. Thus (\ref{JHc}) follows.
	
	\medskip
	\no 
	{\it{Claim}:} There exists a sequence ${(x_n,y_n)}_{n \geq 1}\subseteq K_{H,c}^+ \cap V_R^-$ with $\lvert x_n \rvert \geq \lvert y_n\rvert\geq n$ for all $n\geq 1$.
	
	\medskip
	If no such sequence exists, then we can choose a sequence $(x_n,y_n)\in K_{H,c}^+\cap V_R^-$ such that  $\lvert y_n \rvert$ is bounded by a fixed constant $M>1$ for all $n\geq 1$  and $\lvert x_n \rvert \rightarrow \infty$ as $n\rightarrow \infty$. Without loss of generality, we choose $R>0$ sufficiently large such that $K_{H,c}^+$ and $K_{H,dc}^+$ both are contained in $V_R \cup V_R^-$ where $d=\deg H$. Suppose that the H\'{e}non map is of the form: $H:(x,y)\mapsto (y, p(y)-\delta x)$ with $\delta\neq 0$. Then there exists a subsequence $\{(x_{n_k}, y_{n_k})\}\subset  K_{H,c}^+ \cap V_R^-$  such that $\{H(x_{n_k}, y_{n_k})\}\subset  K_{H,dc}^+ \cap V_R^-$ and thus
	\[
	\lvert y_{n_k} \rvert \geq \lvert p(y_{n_k})-\delta x_{n_k}\rvert  \geq \lvert \delta\rvert\lvert x_{n_k} \rvert-\lvert p(y_{n_k})\rvert.
	\]
Therefore, the sequence $\{x_{n_k}\}$ is bounded  which is a contradiction.
	
	\medskip
	Since $\overline{K_{H,c}^+}=K_{H,c}^+ \cup I^+$ in $\mathbb{P}^2$, 
	\[
	\frac{\vert y_n \vert}{\vert x_n \vert}\rightarrow 0 \text{ as } n\rightarrow \infty. 
	\]
	Thus
	\begin{equation}\label{eps}
	\lvert y_n \rvert \leq \epsilon_n\lvert x_n\rvert
	\end{equation}
	for all $n\geq 1$ with $\epsilon_n \rightarrow 0$.
	
	\medskip
	Note that  $\tau a_2^2(x_n,y_n)=(s_2 y_n+r_2, \alpha_2 x_n + \beta_2 y_n + \delta_2)$ and
	\begin{eqnarray}
	\lvert \alpha_2 x_n + \beta_2 y_n +\delta_2\rvert 
	%&\geq& \lvert \alpha_2\rvert \lvert x_n\rvert-\lvert\beta_2\rvert \lvert y_n\rvert-\lvert \delta_2\rvert \nonumber\\
	&\geq& (\lvert \alpha_2 \rvert - \epsilon_n \lvert \beta_2\rvert)\lvert x_n\rvert-\lvert\delta_2\rvert \nonumber\\
	&\geq& \frac{1}{2}\lvert \alpha_2 \rvert  \lvert x_n \rvert-\lvert \delta_2 \rvert  \nonumber\\
	&\geq& \frac{1}{2} \lvert \alpha_2\rvert \lvert y_n \rvert-\lvert \delta_2 \rvert \geq \lvert s_2\rvert \lvert y_n\rvert + \lvert r_2\rvert \lvert y_n\rvert-\lvert \delta_2\rvert  \nonumber
	\end{eqnarray}
	for all  $n\geq n_0$. Now since  $\lvert s_2\rvert$ and $\lvert r_2\rvert$ can be chosen sufficiently small, we choose them such that ${\lvert\alpha_2\rvert}/{2}\geq\lvert s_2\rvert+ \lvert r_2\rvert$. Thus the last inequality follows.
	
	\medskip 
	Since $\lvert y_n\rvert \rightarrow \infty$ as $n\rightarrow \infty$, it follows that
	\[
	\lvert \alpha_2 x_n + \beta_2 y_n +\delta_2\rvert \geq \lvert s_2 \rvert \lvert y_n\rvert + \lvert r_2\rvert \geq \lvert s_2 y_n + r_2\rvert
	\]
	and 
	\[
	\lvert \alpha_2 x_n + \beta_2 y_n +\delta_2\rvert \geq R
	\]
	for sufficiently large $n$.
	
	\medskip 
	Thus for a sequence ${(x_n,y_n)}_{n\geq 1}\subseteq K_{H,c}^+\cap V_R^-$ with $\lvert x_n \rvert \geq \lvert y_n \rvert \geq n$, it turns out that $\tau a_2^2 (x_n,y_n)\in V_R^+$ for sufficiently large $n$. Thus 
	\[
	(x_n',y_n')= \tau a_1^2 e_1 a_2^1 \tau a_2^2 (x_n,y_n)\in V_R^+
	\]
	and 
	\[
	\lvert b x_n'+c y_n'\rvert\leq (\lvert b\rvert+\lvert c\rvert)\lvert y_n'\rvert.
	\]  
	Hence
	\begin{equation}\label{contra}
	\lvert y_n''\rvert \geq \frac{1}{(\lvert b \rvert+ \lvert c \rvert)} \lvert x_n''\rvert
	\end{equation}
	where $(x_n'',y_n'')=a_1^1(x_n',y_n')$. Now since $F(K_{H,c}^+)=K_{H,c}^+$, 
	\[
	(x_n'',y_n'')=F(x_n,y_n)\in K_{H,c}^+\cap V_R^-
	\]
	for sufficiently large $n\geq 1$ and $\Vert (x_n'',y_n'')\Vert \rightarrow \infty$ as $n\rightarrow \infty$. By (\ref{eps}), we get that 
	\[
	\lvert y_n''\rvert \leq \epsilon_n \lvert x_n''\rvert
	\]
	where $\epsilon_n\rightarrow 0$ as $n\rightarrow \infty$ which contradicts (\ref{contra}). Thus $F$ cannot be of this form.
	
	\medskip
	\no 
	{\it{Case (iii):}} Let 
	\[
	F=e_1\circ a_1 \circ e_2 \circ a_2 \circ \cdots \circ e_k \circ a_k
	\]
	for some $k\geq 1$. Note that $F^{-1}$ has a form as in Case 1. Since $F^{-1}$ also keeps $K_{H,c}^+$ invariant, it follows that $F^{-1}$ is a regular map. Hence, $F$ is a regular map with $[w:1:0]$ as its indeterminacy point for some $w\in \mathbb{C}$.
	
	\medskip
	\no 
	{\it{Case (iv):}} Let 
	\[
	F=e_1\circ a_1 \circ e_2 \circ a_2 \circ \cdots\circ a_{k-1}\circ e_{k-1} \circ e_k 
	\]
	for some $k\geq 1$. For simplicity, we work with 
	\[
	F=e_1\circ a_1 \circ e_2.
	\]
	As in the previous cases, we can write
	\[
	F=e_1 a_1^1\tau a_1^2 e_2
	\]
	and thus,
	\[
	\tau F=\tau e_1 a_1^1\tau a_1^2 e_2.
	\]
	Note that both $\tau e_1 a_1^1$ and $\tau a_1^2 e_2$ are H\'{e}non maps.  Thus $\tau F$ is a H\'{e}non map. Therefore,  $F[0:1:0]=[1:0:0]$ and consequently, $F(V_R^+)$ will intersects $K_{H,c}^+$ since $[1:0:0]$ is the limit point of $K_{H,c}^+$ in $\mathbb{P}^2$. This implies, $V_R^+ \cap K_{H,c}^+\neq \emptyset$ since $F(K_{H,c}^+)=K_{H,c}^+$ . This is clearly a contradiction. Therefore, $F$ can not be of this form.  
	
	\medskip 
	Thus we prove that $F$ is a regular map. Furthermore, the point $[1:0:0]$ is the indeterminacy point either of $F$ or of $F^{-1}$.  Without loss of generality, let $[1:0:0]$ be the indeterminacy point of $F$. 

\begin{center}
\it{$K_{H,c}^+$ never remains  invariant under an automorphism $F$ with $\deg F\geq 2$ and $c>0$:}
\end{center}
From (\ref{Gc1}), it follows that 
\begin{equation*}\label{max}
\max \left \{G_H^+\circ F (x,y)-c,0 \right \}= b^+ \max \left \{G_H^+(x,y)-c,0 \right \}
\end{equation*}
in $\mathbb{C}^2$. Comparing both sides of (\ref{max}), we get that 
\begin{equation*} \label{GFc}
G_H^+ \circ F (x,y)=b^+ G_H^+(x,y)-b^+ c + c
\end{equation*}
for all $(x,y)$
%\in \left \{(x,y):G_H^+(x,y) > c \right \}$, 
such that $G_H^+(x,y)>c$. In particular  for sufficiently large $R>0$, we have that $G_H^+(x,y)>c$ for all  $(x,y)\in V_R^+$.

\medskip 
Since
\begin{equation*} 
G_H^+ \circ F (z)=b^+ G_H^+(z)-b^+ c + c
\end{equation*}
for all $(x,y)$ such that $G_H^+ (x,y)> c$,
%\in \left \{(x,y)\in \mathbb{C}^2:G_H^+ (x,y)> c \right \}$, 
it follows that
$$
F\left(\left \{(x,y)\in \mathbb{C}^2:G_H^+(x,y) > c \right \}\right)\subseteq \left \{(x,y)\in \mathbb{C}^2:G_H^+(x,y) > c \right \}.
$$
Also, 
\begin{equation*} 
G_H^+ \circ F^{-1} (x,y)=b^- G_H^+(x,y)-b^- c + c
\end{equation*}
 for all $(x,y)$
 %\in \left \{(x,y)\in \mathbb{C}^2:G_H^+ (x,y)> c \right \}$ 
 such that $G_H^+(x,y)>c$, it follows that
$$
F^{-1}\left(\left \{(x,y)\in \mathbb{C}^2:G_H^+(x,y) > c \right \}\right)\subseteq \left \{(x,y)\in \mathbb{C}^2:G_H^+(x,y) > c \right \}.
$$
Therefore, 
$$
F\left(\left \{(x,y)\in \mathbb{C}^2:G_H^+(x,y) > c \right \}\right)=\left \{(x,y)\in \mathbb{C}^2:G_H^+(x,y) > c \right \}.
$$
This implies 
$$
F\left(\left \{(x,y)\in \mathbb{C}^2:G_H^+(x,y) < c \right \}\right)=\left \{(x,y)\in \mathbb{C}^2:G_H^+(x,y) < c \right \}
$$
since $F(K_{H,c}^+)=K_{H,c}^+$.

\medskip 
\no 
{\it{Claim:}}
$
\Omega_{H,c}=\left \{(x,y)\in \mathbb{C}^2: G_H^+(x,y) <c \right \} \subseteq K_F^+.
$

\medskip 
Since $F(K_{H,c}^+)=K_{H,c}^+$, by Prop. \ref{relegular attraction}, it follows that $J_{H,c}^+=K_{H,c}^+ \subset K_F^+$. By {\it{rigidity}} results of Dinh--Sibony (see Theorem \ref{DS1} and Theorem \ref{DS2}),  it follows that 
\[
J_{H,c}^+=J_F^+.
\] 
By Prop. \ref{pluri-complex}, the function $G_{H,c}^+$ is the pluricomplex Green function of $J_{H,c}^+$. Again as in Theorem \ref{new rigidity}, we can show that $G_F^+$ is the pluricomplex Green function of $J_F^+$ and  thus 
\[
G_{H,c}^+=G_F^+.
\]
Now since $G_{H,c}^+$ vanishes identically on $\Omega_{H,c}$, it follows that $G_F^+$ also vanishes identically on $\Omega_{H,c}$. Since 
\[
K_F^+=\left\{(x,y): G_F^+(x,y)=0 \right\},
\]
we have that 
\[
\Omega_{H,c}\subseteq K_F^+.
\]
Now by Theorem \ref{DS1}, it follows that for each positive $r<c$, the set $J_{H,r}^+$ supports a positive closed (1,1) current of mass 1. On the other hand $K_F^+$ supports a unique positive closed (1,1) current of mass 1.  Contradiction! This finishes the proof.

\medskip 
Thus $F$ must be an affine automorphism. Since $F(K_{H,c}^+)=K_{H,c}^+$ and $\overline{ K_{H,c}^+}=K_{H,c}^+ \cup I^+$ in $\mathbb{P}^2$,  applying similar arguments as before, we can show that 
$$
F(x,y)=(ax+by+f, dy+g).
$$ 
%with $\lvert a \rvert= \lvert d \rvert=1$.

\section{Proof of the Theorem \ref{two levels}}
Let us first state the following propositon (see \cite{BPK}) which we require to prove Theorem \ref{two levels}. Suppose $H$ is of the form (\ref{henon form}).
\begin{prop}\label{Bot}
	For a given H\'{e}non map $H$, there exist non-vanishing holomorphic functions 
	$\phi_H^\pm: V_R^\pm \rightarrow \mathbb{C}$  such that 
	\[
	\phi_H^+\circ H(x,y)=c_H{(\phi_H(x,y))}^d
	\]
	in $V_R^+$ where
	\[
	c_H = \prod_{j=1}^m {c_j}^{d_{j+1}\cdots d_m}
	\]
	with the convention that $d_{j+1}\cdots d_m=1$ when $j=m$, $d=d_j \cdots  d_1$ and 
	\[
	\phi_H^-\circ H^{-1}(x,y)=c_H'{(\phi_H^-(x,y))}^d
	\]
	in $V_R^-$ where 
	\[
	c_H'=\prod_{j=1}^m {\left({c_j}\delta_j^{-1}\right)}^{d_{j-1}\cdots d_1} 
	\]
	with the convention that $d_{j-1}\cdots d_1=1$ when $j=1$.
	Further,
	\[
	\phi_H^+(x,y)\sim y \text{ as } \lVert(x,y)\rVert\rightarrow \infty \text{ in } V_R^+
	\]
	and
	\[
	\phi_H^-(x,y)\sim x \text{ as } \lVert(x,y)\rVert\rightarrow \infty \text{ in } V_R^-.
	\]
\end{prop}
The functions $\phi_H^\pm$ are called the B\"{o}ttcher functions corresponding to the H\'{e}non map $H$. 
\subsection*{Proof of the Theorem \ref{two levels}:} 
Since $K_{H_1,c_1}^+=K_{H_2,c_2}^+$, by Prop. \ref{pluri-complex}, we have 
$G_{H_1,c_1}^+=G_{H_2,c_2}^+$ in $\mathbb{C}^2$, i.e., 
\begin{equation*}
\max \left \{G_{H_1}^+ -c_1,0\right \}=\max\left \{G_{H_2}^+ -c_2,0\right \}
\end{equation*}
in $\mathbb{C}^2$. Therefore, 
\[
G_{H_1}^+-c_1=G_{H_2}^+-c_2
\]
in $V_R^+$, for $R>0$ sufficiently large.  Thus by Prop. \ref{Bot}, it follows that
\begin{equation} 
\log\left  \lvert \phi_{H_1}^+\right \rvert+\frac{1}{d_{H_1}-1}\log\left  \lvert c_{H_1} \right \rvert=\log \left \lvert \phi_{H_2}^+ \right \rvert+\frac{1}{d_{H_2}-1}\log\left  \lvert c_{H_2} \right\rvert+c
\end{equation}
with $c=c_1-c_2$ in $V_R^+$. Also since 
\begin{equation} \label{phi y}
\phi_{H_i}(x,y )\sim y \text{ as } \lVert (x,y) \rVert \rightarrow \infty
\end{equation}
for $i=1,2$, we have 
\begin{equation}\label{relC}
\frac{1}{d_{H_1}-1}\log\left  \lvert c_{H_1} \right \rvert=\frac{1}{d_{H_2}-1}\log\left  \lvert c_{H_2} \right\rvert+c
\end{equation}
which in turn gives 
\[
\log\left  \lvert \phi_{H_1}^+\right \rvert=\log\left  \lvert \phi_{H_2}^+\right \rvert
\]
in $V_R^+$. Again, using (\ref{phi y}), we have 
\begin{equation}
 \phi_{H_1}^+ = \phi_{H_2}^+=\phi^+ 
\end{equation}
in $V_R^+$.  Using \ref{relC}, we get
\[
\log {\lvert c_{H_1}\rvert}^{d_{H_2}-1}=\log {\lvert c_{H_2}\rvert}^{d_{H_1}-1}+c(d_{H_2}-1)(d_{H_1}-1)
\]
which implies that 
\begin{equation} \label{delta+}
c_{H_1}^{d_{H_2}} c_{H_2}=\delta_+ c_{H_2}^{d_{H_1}} c_{H_1}
\end{equation}
where $\lvert \delta_+ \rvert= e^ {c(d_{H_1}-1) (d_{H_2}-1)}$.

\medskip 
 By Proposition \ref{Bot},  
\begin{equation*}
\phi^+ \circ H_2 \circ H_1 (x,y)=c_{H_2}{\left(\phi^+ \circ H_1 (x,y)\right)}^{d_{H_2}}=c_{H_2} c_{H_1}^{d_{H_2}}{(\phi^+(x,y))}^{d_{H_1} d_{H_2}}
\end{equation*}
and similarly,
\begin{equation*}
\phi^+ \circ H_1 \circ H_2 (x,y)=c_{H_1}{\left(\phi^+ \circ H_2 (x,y)\right)}^{d_{H_1}}=c_{H_1} c_{H_2}^{d_{H_1}}{(\phi^+(x,y))}^{d_{H_1} d_{H_2}}.
\end{equation*}
Thus, 
\begin{equation*}\label{equality}
\phi^+ ( H_2 \circ H_1) =\delta_+ \phi^+ ( H_1 \circ H_2)
\end{equation*}
on $V_R^+$. 
Now since
\begin{equation*}
\phi^+ \circ H_2 \circ H_1 (x,y)\sim {(H_2 \circ H_1)}_2(x,y)  
\end{equation*}
and 
\begin{equation*}
\phi^+ \circ H_1 \circ H_2 (x,y)\sim {(H_1 \circ H_2)}_2(x,y)  
\end{equation*}
as $\lVert(x,y)\rVert\rightarrow \infty$ in $V_R^+$, it follows that for any fix $x_0 \in \mathbb{C}$, 
\[
{(H_2 \circ H_1)}_2(x_0,y)- \delta_+{(H_1\circ H_2 )}_2(x_0,y)\sim 0 \text{ as } \lvert y \rvert \rightarrow \infty.
\]
Since the expression on the left of the above equation is a polynomial in $y$, it follows that 
\[
{(H_2 \circ H_1)}_2(x_0,y)=\delta_+{(H_1 \circ H_2 )}_2(x_0,y)
\]
for all $y\in \mathbb{C}$. Therefore,
\begin{equation}\label{rel1}
{(H_2 \circ H_1)}_2 \equiv \delta_+{(H_1 \circ H_2)}_2
\end{equation}
in $\mathbb{C}^2$. 

\medskip 
 We again use B\"{o}ttcher coordinates to recover the relation between the first components of these maps. As in the  previous case,  similarly one can show that 
\begin{equation*}
\phi_{H_1}^- = \phi_{H_2}^-=\phi^-.
\end{equation*}
Thus using Prop. \ref{Bot}, we get that    
\begin{equation}\label{F1} 
{(c_{H_1}')}^{d_{H_2}} c_{H_2}' {(\phi^- \circ H_1 \circ H_2 (x,y))}^{d_{H_1} d_{H_2}} = \phi^-(x,y)
\end{equation}	
and 
\begin{equation}\label{F2}
{(c_{H_2}')}^{d_{H_1}} c_{H_1}' {(\phi^- \circ H_2 \circ H_1 (x,y))}^{d_{H_1} d_{H_2}} = \phi^-(x,y)
\end{equation}
for all $(x,y)\in {(H_1\circ H_2)}^{-1}(V_R^-)\cap {(H_2\circ H_1)}^{-1}(V_R^-)=\mathcal{A}$. Note that $\cal A$ is an open neighbourhood of $I^+ = [1:0:0]$ in $\mbb P^2$. 
%By (\ref{Green2}), 
%\[
%\vert {(c_H')}^{d_F} c_F'  \vert = \vert {(c_F')}^{d_H} c_H'  \vert
%\]
As in (\ref{delta+}), it can be shown that 
\begin{equation} \label{delta+}
{\left(c_{H_1}'\right)}^{d_{H_2}} c_{H_2}'=\eta {\left(c_{H_2}'\right)}^{d_{H_1}} c_{H_1}'
\end{equation}
where $\lvert \eta\rvert= e^ {d(d_{H_1}-1) (d_{H_2}-1)}$ with $d=d_1-d_2$. Hence
\[
{(\phi^- \circ H_1 \circ H_2 (x,y))}^{d_{H_1} d_{H_2}} = \eta {(\phi^- \circ H_2 \circ H_1 (x,y))}^{d_{H_1} d_{H_2}}
\]
on $\mathcal{A}$. Consequently, there exists $\delta_-$ (an appropriate $d_{H_1} d_{H_2}$-th root of $\eta$) such that
\begin{equation} \label{phi h-}
\phi^-\circ (H_1\circ H_2)=\delta_- \phi^- \circ (H_2 \circ H_1)
\end{equation}
on $\cal A$ where $\vert \delta_- \vert = e^ {d(d_{H_1}-1) (d_{H_2}-1)}$.

\medskip  
Note that for a fixed $c\neq 0$, there exists $\epsilon>0$ sufficiently small such that
\[
\mathcal{A}_{\epsilon,c}=\{[{1}/{y}:c:1] \text{ where }0\neq \lvert y\rvert <\epsilon\}
\]
intersects $\cal A$ and contains $I^+=[1:0:0]$ in its boundary. Choose $[x_n:c:1]\in \cal {A}_{\ep,c}$ such that  $\lvert x_n \rvert \rightarrow \infty$. Now  since $(H_2 \circ H_1)(x_n,c), (H_1\circ H_2)(x_n,c)\in V_R^-$, it follows that 
%$(x_n,c)$ being in $\mathcal{U}$, we have  
\[
{(H_2 \circ H_1)}_1(x_n,c),{(H_1\circ H_2)}_1(x_n,c)\rightarrow \infty
\]
as $n\rightarrow \infty$. 

\medskip
\no 
%Since $\phi^-(x,y)\sim x$ as $\lVert(x,y)\rVert\rightarrow \infty$,  
Using (\ref{phi h-}), 
\begin{equation*}
{(H_1\circ H_2)}_1(x_n,c)-\delta_-{(H_2\circ H_1)}_1(x_n,c)\rightarrow 0
\end{equation*} 
as $n\rightarrow \infty$. The expression on the left is a polynomial in $x$ for each fixed $c$ and thus
\begin{equation*}
{(H_1\circ H_2)}_1(x,c)=\delta_-{(H_2 \circ H_1)}_1(x,c)
\end{equation*}
for all $x\in \mathbb{C}$. Using the similar argument as in the previous case, we get
\begin{equation}\label{rel2}
{(H_2\circ H_1)}_1 \equiv \delta_-{(H_1\circ H_2)}_1
\end{equation}
in $\mathbb{C}^2$.

\medskip 
Hence using (\ref{rel1}) and (\ref{rel2}), we get
\[
H_2 \circ H_1= C\circ H_1\circ H_2
\]
where $C(x, y) = (\delta_- x,\delta_+ y)$ with $\lvert \delta_+\rvert=e^ {c(d_{H_1}-1) (d_{H_2}-1)}$ and  $\lvert \delta_-\rvert=e^ {d(d_{H_1}-1) (d_{H_2}-1)}$ with $c=c_1-c_2$ and $d=d_1-d_2$.

\end{document}